\documentclass[12pt]{amsart}
\usepackage[francais]{babel}

\title{Sur une erreur dans les EGA,\\  \tiny Les \'El\'ements de G\'eom\'etrie Alg\'ebrique} 

\author{Labib Haddad}
\address{120 rue de Charonne, 75011 Paris, France}
\email{labib.haddad@wanadoo.fr}

\usepackage{amssymb}
\usepackage{amsmath}

\newcommand{\su}{\subsection*}

\newcommand{\Ž}{\'e}

\newcommand{\ˆ}{\`a}

\newcommand{\nts}{\negthickspace}

\begin{document}
\maketitle

\markboth{Labib Haddad}{Sur les EGA }

\thispagestyle{empty}

\su{On a small hitch in \lq\lq Les \'El\'ements de G\'eom\'etrie Alg\'ebrique"} 

\

\

\noindent In EGA I, [1], {\sc Grothendieck} and {\sc Dieudonn\Ž} state the following result, on page 59\nts:

\su{Corollary (2.4.4)}\sl Let $X$ be a  N\oe therian space, and ($E_\alpha$) be a directed family of globally constructible subsets of $X$ such that\nts:

$1$. $X$ is the union of the family ($E_\alpha$).

$2$. Each irreducible closed subset of  $X$ is contained in the closure of one of the  $E_\alpha$'s.

There exists, then, an  index $\alpha$ such that $X = E_\alpha$.

Whenever each irreducible closed subset of  $X$ has a generic point, the result still holds with condition 2 deleted.\rm

\

\centerline{\bf This corollary is not correct as the following example shows.}

\su{Counterexample} Let $X = \{1,2,3, \dots,n,\dots, \omega\}$ be the topological space whose set of open subsets is
$$\{\emptyset, \{\omega\}\} \cup \{A : \omega \in A \ \text{and $X\setminus A$ is finite}\}.$$ 
Space $X$ is N\oe therian and each of its singletons $\{x\}$ is  locally closed. Let ($E_\alpha$) be the family of finite subsets of $X$. This is a directed family of globally constructible subsets which cover $X$, such that $\overline{\{\omega\}}= X$. Therefore, conditions 1 and 2 in the  corollary are satisfied, but its conclusion is not.

\

It seems very likely that the error has already been noticed and reported, previously, but I know not when nor where. I apologize to all those whom I could not quote, due to ignorance.

\

I tried, with no avail, to find out if the tear  had spread in the book or not. It would certainly be important to know, but I must leave that to more competent and better equipped than me, than I.

\

\section*{Le texte en fran\c cais}

\

\`A la suite d'une r\Žcurrence n\oe th\Žrienne mal conduite, 
dans EGA I, [1], en page 59, {\sc Grothendieck} et {\sc Dieudonn\Ž} ont pu \Žnoncer le r\Žsultat (faux) suivant.

\su{Corollaire (2.4.4)}\sl Soient $X$ un espace n\oe th\Žrien, ($E_\alpha$) une famille filtrante croissante de parties globalement constructibles de $X$ telle que :

$1$. $X$ est r\Žunion de la famille ($E_\alpha$).

$2$. Toute partie ferm\Že irr\Žductible de $X$ est contenue dans l'adh\Žrence d'un des $E_\alpha$.

Alors il existe un indice $\alpha$ tel que $X = E_\alpha$.

Lorsque toute partie ferm\Že irr\Žductible de $X$ admet un point g\Žn\Žrique, l'hypoth\se 2 peut \^etre supprim\Že.\rm

\

\centerline{\bf Ce corollaire est faux,  le \emph{contre-exemple} suivant le prouve.}

\su{Contre-exemple} Soit $X = \{1,2,3, \dots,n,\dots, \omega\}$. On munit $X$ de la topologie dont l'ensemble des ouverts est 
$$\{\emptyset, \{\omega\}\} \cup \{A : \omega \in A \ \text{et $A$ est cofini dans $X$}\}.$$ 
L'espace $X$ est n\oe th\Žrien et chaque partie $\{x\}$ est localement ferm\Že. On prend pour famille ($E_\alpha$) la famille des parties finies de $X$. C'est une famille filtrante croissante de parties globalement constructibles qui recouvre $X$ et l'on a $\overline{\{\omega\}}= X$. Donc les conditions 1 et 2 du corollaire sont bien remplies, mais la conclusion du corollaire est fausse.

On voudra bien se souvenir  que, dans un espace n\oe th\Žrien, les parties globalement constructibles sont {\it les r\Žunions finies des parties localement ferm\Žes\rm : voir EGA I [1, (2.4.1)].

\su{Remarque} L'erreur vient de ce que la r\Žcurrence n\oe th\Žrienne de la d\Žmonstration est mal conduite. On ne peut pas supposer que toute partie  ferm\Že propre de $X$ est contenue dans un $E_\alpha$ : s'il existe un $Y$ irr\Žductible dont toute partie ferm\Že propre est contenue dans un 
$E_\alpha$, il existe bien $\beta$ tel que l'on ait $Y = \overline{E_\beta}$, mais il se peut tr\s bien que l'on ait $E_\beta \cap Y =\emptyset$ !

\

Sans doute, cette erreur a-t-elle d\Žj\ˆ \Žt\Ž relev\Že (plus d'une fois) et corrig\Že, mais je ne sais pas o\`u ni par qui, ni quand. D'avance, je pr\Žsente toutes mes excuses \ˆ tous ceux que je n'ai pas pu citer, faute de les conna\^itre. 

\

Une vaine recherche ne m'a pas permis de d\Žterminer si cette erreur a pu ou non se propager. Il serait pourtant utile de le savoir et qu'une \Žtude plus pouss\Že soit faite, par des chercheurs plus comp\Žtents, et mieux arm\Žs.

\

{\bf Pour la petite histoire} : un petit nombre de coll\gues Libanais, r\Žunis en s\Žminaire, \ˆ Beyrouth, avaient entrepris d'\Žtudier les EGA, dans les ann\Žes 70, avant le d\Žsastre qui s'est abattu sur le Liban. Par lettre, l'un d'eux m'avait demand\Ž de les aider \ˆ d\Žmontrer ce corollaire. Ces jours derniers, c'est en remuant de nouveau de vieux papiers poussi\Žreux pour les ranger que j'ai retrouv\Ž ma petite note manuscrite, en r\Žponse, et j'ai la faiblesse de vouloir l'afficher sur {\tt arXiv}.

\

\section*{R\Žf\Žrence}

\

[1] {\bf Grothendieck, A.; Dieudonn\Ž, J. A.} \'El\Žments de g\Žom\Žtrie algbrique. I. (French) [Elements of algebraic geometry. I] Grundlehren der Mathematischen Wissenschaften [Fundamental Principles of Mathematical Sciences], 166. Springer-Verlag, Berlin, 1971. ix+466 pp. ISBN\nts: 3-540-05113-9; 0-387-05113-9 {\bf MR3075000}

\

\

\enddocument